%%%%%%%   Geometry and Topology Monographs Volume 1 :
%%%%%%%   m1-16.tex :  Marden for Epschrift : Plain TeX

\input gtmacros
\input gtmonout
\volumenumber{1}
\volumeyear{1998}
\volumename{The Epstein birthday schrift}
\pagenumbers{335}{340}
\received{1 June 1998}
\published{27 October 1998}
\papernumber{16}
\def\ssp{\stdspace}

\title{Complex projective structures on Kleinian groups}

\author{Albert Marden}
\address{School of Mathematics, University of Minnesota\\
Minneapolis, MN 55455, USA}
\email{am@geom.umn.edu}
\abstract
Let $M^3$ be a compact, oriented, irreducible, and boundary
incompressible $3$--manifold.  Assume that its fundamental group is
without rank two abelian subgroups and $\partial M^3 \ne \emptyset$.
We will show that every homomorphism $\theta \co\pi_1(M^3)\to PSL(2,{\bf
C})$ which is not ``boundary elementary" is induced by a possibly
branched complex projective structure on the boundary of a hyperbolic
manifold homeomorphic to $M^3$.
\endabstract
\asciiabstract{%
Let M^3 be a compact, oriented, irreducible, and boundary
incompressible 3-manifold.  Assume that its fundamental group is
without rank two abelian subgroups and its boundary is non-empty.
We will show that every homomorphism from pi_1(M) to PSL(2,C)
which is not `boundary elementary' is induced by a possibly
branched complex projective structure on the boundary of a hyperbolic
manifold homeomorphic to M.}

\primaryclass{30F50}\secondaryclass{30F45, 30F60, 30F99, 30C99}

\keywords{%
Projective structures on Riemann surfaces,
hyperbolic 3--manifolds}
\asciikeywords{%
Projective structures on Riemann surfaces,
hyperbolic 3-manifolds}

\maketitlepage

\section{Introduction}  Let $M^3$ be a compact, oriented, irreducible, and
boundary incompressible $3$--manifold such that its fundamental group
$\pi_1(M^3)$ is without rank two abelian subgroups.  Assume that $\partial
M^3 = R_1\cup\dots\cup R_n$ has $n\geq 1$ components, each  a surface
necessarily of genus exceeding one.

We will study homomorphisms
$$
\theta \co\pi_1(M^3)\to G\subset PSL(2,{\bf C})
$$
onto groups $G$ of M\"obius transformations.  Such a homomorphism is called
{\it elementary} if its image $G$ fixes a point or pair of points in
its action on ${\bf H}^3\cup\partial {\bf H}^3$, ie on hyperbolic $3$--space
and its ``sphere at infinity''.  More particularly, the homomorphism
$\theta$ is called {\it boundary elementary} if the image $\theta
(\pi_1(R_k))$ of some boundary subgroup is an elementary group.  (This
definition is independent of how the inclusion
$\pi_1(R_k)\hookrightarrow\pi_1(M^3)$ is taken as the images of
different inclusions of the same boundary group are conjugate in $G$).

The purpose of this note is to prove:
\proclaim{Theorem 1} Every homomorphism $\theta \co\pi_1(M^3)\to
PSL(2,{\bf C})$ which is not boundary elementary is induced by a possibly
branched complex projective structure on the boundary of some Kleinian
manifold ${\bf H}^3\cup\Omega (\Gamma)/\Gamma\cong M^3$.\endproclaim
This result is based on, and generalizes:
\proclaim{Theorem A}{\rm (Gallo--Kapovich--Marden [1])}\ssp
Let $R$ be a compact, oriented surface of genus exceeding one.  Every
homomorphism $\pi_1(R)\to PSL(2,{\bf C})$ which is not elementary is
induced by a possibly branched complex projective structure on ${\bf
H}^2/\Gamma\cong R$ for some Fuchsian group $\Gamma$. \endproclaim
Theorem 1 is related to Theorem A as simultaneous uniformization is related
to uniformization.  Its application to quasifuchsian manifolds could be
called simultaneous projectivization. For Theorem A finds a single
surface on which the structure is determined whereas
Theorem 1 finds a structure simultaneously on the 
pair of surfaces arising from some quasifuchsian group.
\section{Kleinian groups}  
Thurston's hyperbolization theorem [3] implies that
$M^3$ has a hyperbolic structure: there is a Kleinian group
$\Gamma_0\cong\pi_1(M^3)$ with regular set $\Omega
(\Gamma_0)\subset\partial {\bf H}^3$ such that ${\cal M}(\Gamma_0)={\bf
H}^3\cup\Omega (\Gamma_0)/\Gamma_0$ is homeomorphic to $M^3$.  The group
$\Gamma_0$ is not uniquely determined by $M^3$, rather $M^3$ determines
the deformation space ${\cal D}(\Gamma_0)$  (taking a fixed $\Gamma_0$ as
its origin). 

We define ${\cal D}^*(\Gamma_0)$ as the set of those isomorphisms 
$\phi\co \Gamma_0 \to
\Gamma \subset PSL(2, \bf{C})$ onto Kleinian groups $\Gamma$ which are induced
by orientation preserving homeomorphisms ${\cal M}(\Gamma_0)\to{\cal M}(\Gamma)
$. Then ${\cal D}(\Gamma_0)$ is defined as 
${\cal D}^*(\Gamma_0)/PSL(2,{\bf C})$, since we do not distinguish between
elements of a conjugacy class.

Let ${\cal V}(\Gamma_0)$ denote the representation space
${\cal V}^*(\Gamma_0)/PSL(2,{\bf C})$ where $ {\cal V}^*(\Gamma_0)$ is the
space of boundary nonelementary homomorphisms 
$\theta\co \Gamma_0 \to PSL(2,{\bf C})$.

By Marden [2], ${\cal D}(\Gamma_0)$ is a complex manifold of dimension
$\sum [3({\rm genus}\ R_k)-3]$ and an open subset of the representation
variety ${\cal V}(\Gamma_0)$.
If $M^3$ is acylindrical, ${\cal D}(\Gamma_0)$ is relatively compact in
${\cal V}(\Gamma_0)$ (Thurston [4]).

The fact that ${\cal D}(\Gamma_0)$ is a manifold depends on a uniqueness
theorem (Marden [2]).  Namely two isomorphisms
$\phi_i\co\Gamma_0\to\Gamma_i,\ i=1,2$, are conjugate if and only if
$\phi_2\phi_1^{-1}\co \Gamma_1\to\Gamma_2$ is induced by a homeomorphism
${\cal M}(\Gamma_1)\to {\cal M}(\Gamma_2)$ which is homotopic to a
conformal map.
\section{Complex projective structures}  For the purposes of this note we will
use the following definition (cf [1]).  A {\it complex projective
structure} for the Kleinian group $\Gamma$ is a locally univalent
meromorphic function $f$ on $\Omega (\Gamma)$ with the property that
$$
f(\gamma z) =\theta(\gamma)f(z),\ z\in\Omega (\Gamma),\
\gamma\in\Gamma ,
$$
for some homomorphism $\theta \co\Gamma\to PSL(2,{\bf C})$.  We are free to
replace $f$ by a conjugate $AfA^{-1}$, for example to normalize $f$ on one
component of $\Omega (\Gamma)$.

Such a function $f$ solves a Schwarzian equation
$$
S_f(z) = q(z),\ \ q(\gamma z)\gamma'(z)^2=q(z);\
\gamma\in\Gamma ,\ z\in\Omega(\Gamma), 
$$
where $q(z)$ is the lift to $\Omega (\Gamma)$ of a holomorphic quadratic
differential defined on each component of $\partial {\cal M}(\Gamma)$.
Conversely, solutions of the Schwarzian,
$$
S_g(z) = q(z),\ z\in\Omega (\Gamma) ,
$$
are determined on each component of $\Omega (\Gamma)$ only up to post
composition by any M\"obius transformation.  The function $f$ has the
property that it not only is a solution on each component, but that its
restrictions to the various components fit together to determine a
homomorphism $\Gamma\to PSL(2,{\bf C})$.  Automatically (cf [1]), the
homomorphism $\theta$ induced by $f$ is boundary nonelementary.

When {\it branched} complex projective structures for a Kleinian group
are required, it suffices to work with the simplest ones: $f(z)$ is
meromorphic on $\Omega (\Gamma)$, induces a homomorphism $\theta
\co\Gamma\to PSL(2,{\bf C})$ (which is automatically boundary
nonelementary), and is locally univalent except at most for one point,
modulo Stab$(\Omega_0)$, on each component $\Omega_0$ of
$\Omega(\Gamma)$.  At an exceptional point, say $z=0$,
$$
f(z) = \alpha z^2(1+o(z)),\ \alpha\ne 0 .
$$
Such $f$ are characterized by Schwarzians with local behavior
$$
S_f(z) = q(z) = -3/2z^2+b/z+\Sigma a_iz^i,\qquad  b^2+2a_0 = 0 .
$$

At any designated point on a component $R_k$ of $\partial {\cal
M}(\Gamma)$, there is a quadratic differential with leading term
$-3/2z^2$.  To be admissible, a differential must be the sum of this and
any element of the $(3g_k-2)$--dimensional space of quadratic differentials
with at most a simple pole at the designated point.  In addition it must
satisfy the relation $b^2+2a_0 = 0$.  That is, the admissible differentials
are parametrized by an algebraic variety of dimension $3g_k-3$.
For details, see [1].

If a branch point needs to be introduced on a component $R_k$ of $\partial
{\cal M}(\Gamma)$, it is done during a construction.  According to [1], a
branch point needs to be introduced if and only if the restriction
$$
\theta \co\pi_1(R_k)\to PSL(2,{\bf C})
$$
does {\it not} lift to a homomorphism
$$
\theta^*\co\pi_1(R_k)\to SL(2,{\bf C}) .
$$
\section{Dimension count}  
The vector bundle of holomorphic quadratic differentials over the
Teichm\"uller space of the component $R_k$ of $\partial {\cal M}(\Gamma_0)$
has dimension $6g_k-6$.  All together these form the vector bundle ${\cal
Q}(\Gamma_0)$ of quadratic differentials over the Kleinian deformation
space ${\cal D}(\Gamma_0)$.  That is, ${\cal Q}(\Gamma_0)$ has
{\it twice} the dimension of ${\cal V}(\Gamma_0)$. The count remains the
same if there is a branching at a designated point.

For example, if $\Gamma_0$ is a quasifuchsian group of genus $g$, ${\cal
Q}(\Gamma_0)$ has dimension $12g-12$ whereas ${\cal V}(\Gamma_0)$ has
dimension $6g-6$.  Corresponding to each non-elementary homomorphism
$\theta \co\Gamma_0\to PSL(2 ,{\bf C})$ that lifts to $SL(2,{\bf C})$ is a
group $\Gamma$ in ${\cal D}(\Gamma_0)$ and a quadratic differential on the
designated component of $\Omega(\Gamma).$ This in turn determines a differential
on the other component.  There is a
solution of the associated Schwarzian equation $S_g(z) = q(z)$ satisfying
$$
f(\gamma z) = \theta (\gamma)f(z),\ z\in\Omega (\Gamma),\
\gamma\in\Gamma .
$$

Theorem 1 implies that ${\cal V}(\Gamma_0)$ has at most $2^n$ components.  
For this is the number of
combinations of $(+,-)$ that can be assigned to the $n$--components of
$\partial {\cal M}(\Gamma_0)$ representing whether or not a given
homomorphism lifts. For a quasifuchsian group $\Gamma_0$,  
${\cal V}(\Gamma_0)$ has two components (see [1]).

\section{Proof of Theorem 1}  We will describe how the construction introduced in
[1] also serves in the more general setting here.

By hypothesis, each component $\Omega_k$ of $\Omega (\Gamma_0)$ is simply
connected and covers a component $R_k$ of $\partial {\cal M}(\Gamma_0)$.
In addition, the restriction
$$
\theta \co\pi_1(R_k)\cong Stab(\Omega_k)\to G_k\subset PSL(2,{\bf C})
$$
is a homomorphism to the nonelementary group $G_k$.

The construction of [1] yields a simply connected Riemann surface ${\cal J}_k$
lying over $S^2$,
called a pants configuration, such that:

(i)\ssp There is a conformal group $\Gamma_k$ acting freely in ${\cal J}_k$
such that ${\cal J}_k/\Gamma_k$ is homeomorphic to $R_k$.

(ii)\ssp  The holomorphic projection $\pi \co{\cal J}_k\to S^2$ is locally
univalent if $\theta$ lifts to a homomorphism $\theta^*\co\pi_1(R_k)\to
SL(2,{\bf C})$.  Otherwise $\pi$ is locally univalent except for one branch
point of order two, modulo $\Gamma_k$.

(iii)\ssp  There is a quasiconformal map $h_k\co\Omega_k\to {\cal J}_k$ such that
$$
\pi h_k(\gamma_z) = \theta (\gamma)\pi h_k(z),\ \gamma\in
Stab(\Omega_k),\ z\in\Omega_k .
$$

Once $h_k$ is determined for a representative $\Omega_k$ for each component
$R_k$ of $\partial {\cal M}(\Gamma_0)$, we bring in the action of
$\Gamma_0$ on the components of $\Omega (\Gamma_0)$ and the corresponding
action of $\theta (\Gamma_0)$ on the range.  By means of this action a
quasiconformal map $h$ is determined on all $\Omega (\Gamma_0)$ which
satisfies
$$
\pi h(\gamma z) =\theta (\gamma)\pi h(z),\ \gamma\in\Gamma_0,\
z\in\Omega (\Gamma_0) .
$$

The Beltrami differential $\mu (z) = (\pi h)_{\bar z}/(\pi h)_z$ satisfies
$$
\mu (\gamma z)\bar\gamma'(z)/\gamma'(z) = \mu (z),\ \gamma\in\Gamma_0,\
z\in\Omega (\Gamma_0) .
$$
It may equally be regarded as a form on $\partial {\cal M}(\Gamma_0)$.
Using the fact that the limit set of $\Gamma_0$ has zero area, we can
solve the Beltrami equation $g_{\bar z} = \mu g_z$ on $S^2$.  It has
a solution which is a quasiconformal mapping $g$ and is 
uniquely determined up to post composition with a M\"obius
transformation. Furthermore $g$ uniquely determines, up to conjugacy,
 an isomorphism $\varphi
\co\Gamma_0\to\Gamma$ to a group $\Gamma$ in ${\cal D}(\Gamma_0)$.

The composition $\pi h g^{-1}$ is a meromorphic function on each
component of $\Omega (\Gamma)$.  It satisfies
$$
(\pi h g^{-1})(\gamma z) = \theta\varphi^{-1}(\gamma)\pi
h g^{-1}(z),\ \ \gamma\in\Gamma ,\ z\in\Omega (\Gamma) .
$$
The composition is locally univalent except for at most one point on each
component of $\Omega (\Gamma)$, modulo its stabilizer in $\Gamma$.  That
is, $\pi\circ h\circ g^{-1}$ is a complex projective structure on $\Gamma$
that induces the given homomorphism $\theta$, via the identification
$\varphi$.
\section{Open questions}  Presumably, a nonelementary homomorphism $\theta
\co\Gamma_0\to PSL(2,{\bf C})$ can be elementary for one, or all, of the
$n\geq 1$ components of $\partial {\cal M}(\Gamma_0)$.  Presumably too,
the restrictions to $\partial {\cal M}(\Gamma_0)$ of a boundary
nonelementary homomorphism can lift to a homomorphism into $SL(2,{\bf C})$
without the homomorphism $\Gamma_0\to PSL (2,{\bf C})$ itself lifting.
However we have no examples of these phenomena.

According to Theorem 1, there is a subset ${\cal P}(\Gamma_0)$ of
the vector bundle ${\cal Q}(\Gamma_0)$ consisting of those
homomorphic differentials giving rise to, say, unbranched complex projective
structures on the groups in ${\cal D}(\Gamma_0)$.
What is the analytic structure of ${\cal P}(\Gamma_0)$; is it a 
nonsingular, properly embedded, analytic subvariety?

When does a given
Schwarzian equation $S_f(z) = q(z)$ on $\Omega (\Gamma)$ have a solution
which induces a homomorphism of $\Gamma$?  

\penalty-800\vskip-\lastskip\vskip 15pt plus10pt minus5pt

{{\large\bf References}\ppar      % `References' is set \large bold 
\leftskip=25pt\frenchspacing      % The list of references is set 
\parskip=3pt plus2pt\small        % \small  with small spaces between,
\def\ref#1#2\par{\noindent        % numbers in [,]'s and set just to the
\llap{[#1]\stdspace}#2\par}       % left of a 25pt margin.
\def\,{\thinspace}

\ref{1}{\bf D Gallo}, {\bf M Kapovich}, {\bf A Marden}, {\it The
monodromy groups of Schwarzian equations on compact Riemann surfaces},
preprint (1997 revised)  

\ref{2}{\bf A Marden}, {\it The geometry of finitely generated
Kleinian groups}, Annals of Math. 99 (1974) 383--462

\ref{3}{\bf W\,P Thurston}, {\it Three dimensional manifolds,
Kleinian groups, and hyperbolic geometry}, Bull. Amer. Math. Soc. 6 
(1982) 357--381

\ref{4}{\bf W\,P Thurston}, {\it Hyperbolic structures on
$3$--manifolds I}, Annals of Math. 124 (1986) 203--246

}

\Addresses\recd

\end